\documentclass[final,5p,times]{elsarticle}

\usepackage{graphics} 
\usepackage{epsfig} 
\usepackage{mathptmx} 
\usepackage{amsmath} 
\usepackage{amssymb}  
\usepackage{mathtools}
\usepackage{caption}
\usepackage{subcaption}
\usepackage{multirow}
\newcommand{\ub}{\mathbf{u}}
\newcommand{\yb}{\mathbf{y}}

\usepackage{color}
\usepackage{url}
\usepackage{hyperref}

\newtheorem{theorem}{Theorem}
\newtheorem{remark}{Remark}

\newtheorem{lemma}{Lemma}
\newproof{pf}{Proof}
\newproof{pot}{Proof of Theorem \ref{thm}}
\newcommand{\red}{\textcolor{black}}

\journal{Systems and Control Letters}

\begin{document}

\begin{frontmatter}
\title{Exact representation and efficient approximations of linear model predictive control laws via HardTanh type deep neural networks }

\author[upb]{Daniela Lupu}

\affiliation[upb]{organization={Department of Automatic Control and Systems Engineering, University Politehnica Bucharest},
            addressline={}, 
            city={Bucharest},
            postcode={060042}, 
            state={},
            country={Romania}}
            
\author[upb,imsam]{ Ion Necoara}
\affiliation[imsam]{ organization ={Gheorghe Mihoc-Caius Iacob Institute of Mathematical Statistics and Applied Mathematics of the Romanian Academy},
city={Bucharest},
postcode={050711},
country={Romania}}

\begin{abstract}
Deep neural networks have revolutionized many fields, including image processing, inverse problems, text mining  and more recently, give very promising results in systems and control. Neural networks with hidden layers  have a strong potential as an approximation framework of predictive control laws as they usually yield better approximation quality and  smaller memory requirements than existing explicit (multi-parametric) approaches. In this paper, we first show that  neural networks with HardTanh activation functions can exactly represent  predictive control laws of linear time-invariant systems. We derive theoretical bounds on the minimum number of hidden layers and neurons that a  HardTanh neural network should have to exactly represent a given predictive control law. The choice of HardTanh deep neural networks is particularly suited for linear predictive control laws as they usually require less hidden layers and neurons than deep neural networks with ReLU units for representing exactly continuous piecewise affine (or equivalently min-max) maps. In the second part of the paper we bring the physics of the model and standard optimization techniques into the architecture design, in order to eliminate the disadvantages of the black-box HardTanh learning. More specifically,  we  design trainable unfolded HardTanh deep architectures for learning linear predictive control laws  based on two standard iterative optimization algorithms, i.e., projected gradient descent  and accelerated projected gradient descent.  We also study the performance  of the proposed HardTanh type deep neural networks on a  linear model predictive control  application.
\end{abstract}


\end{frontmatter}


\section{Introduction}

\noindent Model predictive control (MPC) is an advanced control
technology that is implemented for a variety of systems and is popular due to its ability to handle  input and state constraints \cite{RawMay:17}. MPC is based on the repeated solution of constrained optimal control problems.  Due to its constraint handling capabilities,
MPC can outperform classic control strategies in many  applications. 

\medskip 

\noindent In its implicit variant, the optimal control inputs are computed by solving an optimization problem on-line that minimizes some cost function subject to state/input constraints, taking the current state as an initial condition. In particular, for linear MPC, the corresponding optimal control problem
can be recast as a convex quadratic program (QP).  Many numerical algorithms
have been  developed in the literature that exploit efficiently the
structure arising in QPs. Basically,  we can
identify three popular classes of optimization algorithms to solve QPs: active set methods \cite{BarBie:06,FerKir:14}, interior point methods \cite{MatBoy:09,DomZgr:12} and (dual) first order methods \cite{LanMon:08,NecNed:13,NedNec:12,PatBem:12,StaSzu:14}.

\medskip 

\noindent The optimal solutions to these MPC-based QPs can also be pre-computed off-line using multi-parametric quadratic programming \cite{BemMor:02,mpt},  with the initial condition of the system viewed as a parameter in the QPs. Therefore, one can also implement an explicit version of the MPC policy,  which consists of implementing a piecewise affine (PWA) control law. Then, in the on-line phase, the main computational effort consists of the identification over some polyhedral partition of the state-space of the region in which the system state resides, e.g., using a binary search tree, and implementing the associated control action.

\medskip 

\noindent \textit{Related work:} Although the theory underlying MPC is quite mature,  implementations of both implicit and explicit MPC suffer from computational difficulties.  In the case of implicit MPC, the computational effort of solving the QPs in real-time makes it impractical for fast systems. Conversely, in explicit MPC  the complexity of the associated PWA control law, measured by the number of affine pieces and regions, grows exponentially with the  number of states of the system and  constraints \cite{BemMor:02}, making it intractable for large systems. Therefore, as an alternative one can consider approximating the exact implicit/explicit MPC policy. Deep neural networks (DNNs) are particularly attractive in approximating the linear  MPC law (i.e., a PWA/min-max map), in view of their universal approximation capabilities, typically requiring a relatively small number of parameters \cite{HorSti:89}. Despite their computationally demanding off-line training requirements, the on-line evaluation of DNN-based approximations to MPC laws is computationally cheap, since it only requires the evaluation of an input-output map \cite{GooBen:16}. Many papers in the MPC  literature  have considered generic deep (black-box) neural networks based on the rectified linear units (ReLU)  activation functions for approximating the MPC policy, see e.g., \cite{MadMor:20,KatPap:20,FG:2021,KarLuc:20}. In particular, in \cite{KarLuc:20} it was shown that ReLU networks can encode exactly the PWA map resulting from the formulation of a linear MPC \red{with politopic  constraints on inputs and states}, with theoretical bounds on the number of hidden layers and neurons required for such an exact representation. However, these results  show that the number of neurons needed in a ReLU type DNN to exactly represent a given PWA function grows exponentially with the number of pieces \cite{chen:22}. Hence, these  papers seem to indicate that the cost of representing a linear MPC law  with  ReLU type DNN is expensive.  

\medskip

\noindent \textit{Contributions:} In this paper, we first show that deep neural networks with HardTanh activation functions can exactly represent  predictive control laws of linear time-invariant systems. We derive theoretical bounds on the minimum number of hidden layers and neurons that a  HardTanh deep neural network should have to exactly represent a given PWA/min-max map, or equivalently a linear predictive control law. {More specifically, we prove that the complexity of the network  is linear in the dimension of the output of the network,  in the number of affine terms and also in the number of min terms. In particular,  we show that the number of hidden neurons required to exactly represent a PWA/min-max function is at most an expression given by the number of pieces times number of min terms.} Therefore,   the choice of HardTanh deep neural networks is particularly suited for linear predictive control laws as they usually require substantially less hidden layers and neurons than ReLU deep neural networks  for representing exactly PWA/min-max maps.

\medskip

\noindent  Although generic HardTanh deep neural networks give promising results in MPC, they are black-box type approaches. Recently unrolled deep learning techniques were proposed, which bring the physics of the system model and standard optimization techniques into the architecture design, in order to eliminate the disadvantages of the deep learning black-box approaches \cite{LePus:22}, \cite{KisOgu:22}. More precisely, the unfolded deep architectures rely on an objective function to minimize, simple constraints,  and an optimization-based algorithmic procedure adapted to the properties of the objective function and constraints. In order to understand  their mechanisms in the context of MPC,  we study the network designed from the unrolled primal first order optimization iterations. To the best of our knowledge,  unfolded architectures for linear MPC and the impact of accelerated schemes on learning performance of MPC policy has not been studied yet and is the core of our second part of the paper. To better understand the mechanisms with or without employing an accelerated scheme and  to facilitate the presentation of the proposed unfolded networks, we focus on two standard  optimization algorithms of the MPC literature: projected gradient descent (PGD) and accelerated projected gradient descent (APGD). \red{Our unfolded HardTanh deep neural networks are expressed as the combination of weights, biases and HardTanh activation functions with expressions depending on the structure of the MPC objective function and constraints, and an extra parameter allowing to consider the basic algorithm (PGD) or its accelerated version (APGD).} The exact/approximate representations of MPC laws  from this paper based on  the generic DNN ($\ub^{ \texttt{\scriptsize HTNN}}$), the unfolded DNN ($\ub^{ \texttt{\scriptsize U-HTNN}}$) and the optimization algorithms  ($\ub^{\texttt{\scriptsize FOM}}$) are summarized in Table 1. Finally, the performance and stability of these representations are illustrate via simulations. 



\begin{table}[]
    \centering
    \begin{tabular}{|l|c|c|}
    \hline
         MPC laws  & exact/approximation & sec. \\
        \hline
         $\ub^{\texttt{\scriptsize FOM}}$ & $ \| \ub^{*}(x_0) - \ub_K^{\texttt{\scriptsize FOM}}(x_0) \| \!\leq\! \mathcal{O} \left(\!(1 \!-\! \kappa^{\texttt{\scriptsize FOM}})^K \!\right)$ &  2 \\
         $\ub^{ \texttt{\scriptsize HTNN}}$ & $\ub^{ \texttt{\scriptsize HTNN}}(x_0) = \ub^*(x_0)$  &  3\\
         $\ub^{\texttt{\scriptsize U-HTNN}}$ & $\| \mathbf{u}^{*}(x_0) -  \ub^{ \texttt{\scriptsize U-HTNN}}(x_0) \| \leq \mathcal{O} \left(\zeta,w \right)$ &  4 \\
        \hline  
         
    \end{tabular}
    \caption{Exact/approximate representations of MPC laws  from this paper: optimization algorithms,  generic DNN and  unfolded DNN. }
    \label{tab:my_label}
\end{table}

\medskip

\noindent  \textit{Contents}. The remainder of this paper is organized as follows. Section
2 provides a brief recall of linear MPC optimization formulation, Section 3 presents explicit
bounds for HardTanh type DNN to be able to represent exactly an MPC law and serves as a motivation for the approximation of MPC laws using the unfolding approach presented in  Section 4.   Section 5 is dedicated to the numerical experiments.


\section{Linear Model Predictive Control (MPC)}
\noindent We consider discrete-time linear systems, defined by the following linear difference equations:
\begin{equation} \label{eq:sysD}
    x_{t+1} = Ax_t + B u_t,
\end{equation}
where $x_t \in \mathbb{R}^{n_x}$ the state at time $t$ and  $u_t \in \mathbb{R}^{n_u}$ the associated control input. We also impose  input constraints:  $u_t \in \mathbb{U}$ for all $t\geq 0$. We assume that $\mathbb{U}$ is a simple polyhedral set, i.e., it is described by linear inequality constraints (e.g., $\underline{u} \leq u \leq \overline{u}$). For system \eqref{eq:sysD}, we consider convex quadratic stage and final costs:
\[ 
\ell(x, u)= 1/2 \left( x^{\top} Q x+u^{\top} R u \right), \quad \ell_f (x)=1/2 x^{\top} P x,  
\] 
respectively. The matrices  $P, Q \in \mathbb{R}^{n_x \times n_x}$  and $R \in \mathbb{R}^{n_u \times n_u}$ are positive (semi-)definite. For a prediction horizon of length  $N$, 
the   MPC problem for system  \eqref{eq:sysD}, given the initial state $x_0$, is:
\
\begin{align} \label{eq:mpc}
&\min _{x_t, u_t} \sum_{t=0}^{N-1} \ell\left(x_t, u_t\right) + \ell_f\left(x_N\right)\\
&\text{s.t.:}  \quad x_{t+1}=A x_t+B u_t, \;\;\;    u_t \in \mathbb{U} \;\;\;  \forall t=0:N-1,  \quad x_0 \; \text{given}. \nonumber 
\end{align}

\noindent In order to present a compact reformulation of the MPC
problem as an optimization problem, we define the 
matrices:
\begin{align} \label{eq:mpc1}
&\mathbf{A} \!=\! \left[\begin{array}{cccc}
B & 0  &  \ldots & 0 \\
A B & B  &  \ldots & 0 \\
A^2 B & A B  &  \ldots & 0 \\
\vdots & \vdots   & \vdots & \vdots \\
A^{N-1} B & A^{N-2} B   & \ldots & B
\end{array}\right],
&& \!\!\!\! \mathbf{A}_N \!=\!\left[\begin{array}{c}
A \\
A^2 \\
A^3 \\
\vdots \\
A^N
\end{array}\right]\!, 
\end{align}
\begin{align*}
& \mathbf{Q} = \text{diag}(\underbrace{Q, \cdots, Q}_{N-1\; \text{times}}, P) \quad \text{and }&& \mathbf{R} = \text{diag}(\underbrace{R, \cdots, R}_{N \; \text{times}}). 
\end{align*}

\noindent Further, by eliminating the state variables, problem \eqref{eq:mpc} is reformulated as a Quadratic Program (QP)  \cite{RawMay:17}:
\begin{align}\label{eq:qp}
    & \min_{\ub \in \mathbb{U}^N } \frac{1}{2} \ub^T \mathbf{H} \ub + (\mathbf{q}(x_0))^T \ub,
\end{align}
where the decision variable contains the entire sequence of commands over the prediction horizon $$\ub = [u_0^T, u_1^T, \cdots, u_{N-1}^T ]^T \in \mathbb{R}^{Nn_u}$$ and 
\[ \mathbf{H} = \mathbf{R} + \mathbf{A}^T \mathbf{Q} \mathbf{A} \quad  \text{and} \quad  \mathbf{q}(x_0) = \mathbf{A}^T \mathbf{Q} \mathbf{A}_N x_0. \]  

\noindent Throughout the paper we assume that the matrices $Q, P$ and $R$ are chosen such that problem \eqref{eq:qp} is strongly convex, i.e.,  $\mathbf{H} \succ 0$. Throughout the paper, we denote  $Nn_u = \nu$. In the next sections we provide several strategies to solve the QP problem \eqref{eq:qp}.


\subsection{Multi-parametric quadratic methods}
\noindent One can easily see that the linear term in the objective function of problem \eqref{eq:qp}, $\mathbf{q}(x_0)$, depends linearly on the initial state of the system $x_0$. Hence, problem \eqref{eq:qp}  can be  viewed as a multi-parametric QP over  a given polytope   $\mathbb{X}$ that defines a certain set of parameters $x_0$ for which a feasible solution exists. If we denote an optimal solution of this (multi-parametric) QP, for a fixed $x_0$,  with $\mathbf{u}^*(x_0)$, then it is known that such a solution, viewed as a function of $x_0$, has a continuous piecewise affine (PWA) representation \cite{BemMor:02}:
\begin{align*}
\mathbf{u}^*(x_0) = \mathcal{P}_{WA}(x_0;\mathbf{H},\mathbf{q},\mathbb{U}) = \begin{cases}
    & G_1 x_0   + g_1, \quad \text{if} \;\; x_0 \in \mathcal{R}_1\\
    & \vdots\\
    & G_m x_0   + g_m, \quad\text{if} \;\; x_0 \in \mathcal{R}_m
\end{cases},    
\end{align*}
\red{where $m$ is the number of regions in  $\mathbb{X}$, i.e., $\mathcal{R}_i$'s are polyhedral sets covering $\mathbb{X}$ ($\mathbb{X} = \cup_{i=1}^m \mathcal{R}_i$).  Note that usually the number $m$ of regions of the piecewise affine solution is, in general, larger than the number of feedback gains, because non-convex critical regions are split into several convex sets \cite{BemMor:02}. However, we can always  associate to each region $\mathcal{R}_i$ an affine feedback gain $G_i x_0   + g_i$ (which may be the same for several regions)}.    Obviously, this PWA map, $ \mathcal{P}_{WA}$, depends on the problem data   $(\mathbf{H},\mathbf{q},\mathbb{U})$. Note that when  the set $\mathbb{U}$ is a box, i.e., $\mathbb{U} =[\underline{u}, \overline{u}] $, the constraints in \eqref{eq:qp} are also  box type:
$$ \underline{\ub} \leq \ub \leq \overline{\ub}, $$
and then  the PWA map depends on the problem data as follows:
\begin{align}
\label{pwl}
\mathbf{u}^*(x_0) = \mathcal{P}_{WA}(x_0;\mathbf{H},\mathbf{q},\underline{\ub}, \overline{\ub}).    
\end{align}
It is known that any continuous piece-wise affine map can be written equivalently as a min-max map \cite{Sch:12} (Proposition 2.2.2):

\begin{lemma}
If $\mathcal{P}_{WA}: \mathbb{X} \rightarrow \mathbb{R}$ is a piecewise affine with affine selection functions $\mathcal{A}_1(x)=c_1^T x+d_1, \ldots, \mathcal{A}_m(x)=c_m^T x+d_m$, then there exists a finite number of index sets $\mathcal{C}_1, \ldots, \mathcal{C}_l \subseteq\{1, \ldots, m\}$ such that
$$
\mathcal{P}_{WA}(x)= \min_{1 \leq i \leq l}  \max_{j \in \mathcal{C}_i} c_j^T x+d_j.
$$
\end{lemma}

\noindent It is shown in \cite{XuBoo:16} that  the number of parameters needed to describe a continuous PWA function is largely reduced by using the min-max representation.  Besides, the online evaluation speed of the function is also improved.   Therefore, in the sequel we consider that the linear MPC law $\mathbf{u}^*(x_0)$ can be  represented as a min-max map (see e.g., \cite{XuBoo:16} for efficient numerical  methods to transform a PWA function into a min-max function): 
\begin{align}
\label{pwl-mm}
\mathbf{u}^*(x_0) &= \mathcal{M}_{in}\mathcal{M}_{ax}(x_0;\mathbf{H},\mathbf{q},\underline{\ub}, \overline{\ub}).
\end{align}
 Note that the piece-wise affine/min-max map, $\mathbf{u}^*(\cdot)$, defined on  a given polytope   $\mathbb{X}$ can be computed off-line, in closed form using multi-parametric QP tools,  such as MPT toolbox \cite{mpt}. Then, on-line, for a given initial state $x_0$, we  need to search in a look-up table to generate $\ub^*(x_0)$. This two phases strategy is called explicit MPC \cite{BemMor:02}.  However, the complexity of computing this piece-wise affine/min-max map, off-line, grows exponentially with the problem data in \eqref{eq:qp}. Also, the search in the on-line phase can be expensive as well. Due to the high complexity of computing explicitly this map $\ub^*(\cdot)$  and of the on-line search, one can consider as an alternative the computation of $\ub^*(x_0)$ for a given $x_0$ on-line using first order methods, which have tight certificates of convergence. In the next section we discuss such on-line methods.


\subsection{First order optimization methods}
\noindent Since $\mathbb{U}$ is a simple polyhedral set, we can consider primal first order methods (\red{FOM}) to solve problem \eqref{eq:qp}, since the projection onto $\mathbb{U}$ is easy (see \cite{Nes:04}). A standard algorithm is the Projected Gradient Descent (PGD) \cite{Nes:04}. Providing an initial point $\ub_0$, the PGD iteration has the expression: 
\begin{align}
\label{pgd}
   &\ub_{k+1} = \Pi_{\mathbb{U}^N}\left( (\operatorname{I}_{\nu} - \alpha_k \mathbf{H}) \ub_{k} - \alpha_k  \mathbf{q}(x_0)\right),
\end{align} 
where $\Pi_{\mathbb{U}^N}$ denotes the  projection operator onto the set $\mathbb{U}^N$,  $\alpha_k$ is an appropriate stepsize (e.g., $\alpha_k = 1/\|\mathbf{H}\|$) and $\operatorname{I}_{\nu}$ denotes the identity matrix of dimension $\nu$. Having $\mathbf{H} \succ 0$, PGD converges linearly with a rate depending on the condition number of the matrix $\mathbf{H}$ \cite{Nes:04}. For faster convergence, one can use the accelerated version of PGD, namely the Accelerated Projected Gradient Descent (APGD) \cite{Nes:04}. The improvement in convergence is payed by growing a little the complexity per iteration by adding an extrapolation step: 
\begin{align}
\label{apgd}
    &\ub_{k+1} = \Pi_{\mathbb{U}^N}\left( (\operatorname{I}_{\nu} - \alpha_k \mathbf{H}) \yb_{k} - \alpha_k  \mathbf{q}(x_0)\right) \\
    &\yb_{k+1} = (1+\beta_k) \ub_{k+1} - \beta_k \ub_{k}, \nonumber 
\end{align}
where $\ub_0 = \yb_0$, $\alpha_k$ is an appropriate stepsize (usually chosen as in PGD) and the extrapolation parameter $\beta_k$ can be chosen even constant when the objective function is strongly convex (i.e., $\mathbf{H} \succ 0$ in \eqref{eq:qp}), see \cite{Nes:04} for more details. One can immediately notice that  if we choose $\beta_k=0$, then we recover PGD algorithm.  It is known that APGD converges faster than PGD, i.e.,  it has a linear  rate depending on the square root  of the condition number of the matrix $\mathbf{H} \succ 0 $ \cite{Nes:04}.  Note that for a box set, i.e., $\mathbb{U} = \left[ \underline{u}, \overline{u}\right]$, the projection operator has the following simple expression:
\begin{align*}
\left[\Pi_{\mathbb{U}^N}(\ub)\right]_i = \begin{cases}
 & \underline{\ub}_i,\quad  \ub_i \leq \underline{\ub}_i \\
& \ub_i,  \quad \underline{\ub}_i< \ub_i < \overline{\ub}_i \\
& \overline{\ub}_i, \quad \ub_i \geq \overline{\ub}_i.
\end{cases}.
\end{align*}

\noindent It is obvious that  min-max functions are closed under addition, multiplication, minimization and maximization operations. \; Hence, combining the min-max expression of the projection operator, $\left[\Pi_{\mathbb{U}^N}(\ub)\right]_i = \min(\overline{\ub}_i, \max(\ub_i, \underline{\ub}_i))$,  with the update rules for the first order methods presented before (PGD and APGD), one can  notice that after $K$ iterations of one of these methods, \eqref{pgd} or \eqref{apgd}, rolling backward, we obtain again a min-max representation for the last iterate $\ub_K$, denoted $\mathbf{u}^{\texttt{\scriptsize FOM}}_K(x_0)$, which depends on the starting point $\ub_0$, initial state $x_0$ and problem data:
\begin{align}
\label{mm-fom}
\mathbf{u}^{\texttt{\scriptsize FOM}}_K(x_0) &= \mathcal{M}_{in}^{\texttt{\scriptsize FOM}}\mathcal{M}_{ax}^{\texttt{\scriptsize FOM}}(x_0;\mathbf{H},\mathbf{q},\underline{\ub}, \overline{\ub},\ub_0,K).    
\end{align}
From the exiting convergence guarantees for first order methods (i.e., for PGD and APGD), one can easily obtain a certificate about how far is the approximate min-max map \eqref{mm-fom} from the optimal min-max map \eqref{pwl-mm}, see \cite{Nes:04}: 
\begin{align}
\label{conv-fom}
\| \mathbf{u}^{*}(x_0) - \mathbf{u}^{\texttt{\scriptsize FOM}}_K(x_0) \|^2 \leq \mathcal{O} \left((1 - \kappa^{\texttt{\scriptsize FOM}})^K \right) \quad \forall x_0 \in \mathbb{X},
\end{align} 
where $\kappa^{\texttt{\scriptsize FOM}} = \kappa$  for PGD and $\kappa^{\texttt{\scriptsize FOM}} = \sqrt{\kappa}$ for APGD (here, $\kappa$ denotes the condition number of the matrix $\mathbf{H}$). Note that the bound \eqref{conv-fom} is usually tight and is uniform in $x_0 \in \mathbb{X}$ (the initial state of the system), provided that $\mathbb{U}$ is bounded set (i.e., in $\mathcal{O}(\cdot)$ there is  no dependence~on~$x_0$, but on the diameter of $\mathbb{U}$).



\section{HardTanh deep neural network (HTNN) for exact representation of  MPC laws}

\noindent In this section, we design the  \texttt{HardTanh} deep neural network (\texttt{HTNN}) arhitecture and show that any given continuous min-max affine function can be exactly represented by a \texttt{HTNN}. We also provide the complexity of such network. Let us start by summarizing some fundamental concepts of deep neural networks. We define a feed-forward neural network as a sequence of layers of neurons which determines a function $f: \mathbb{R}^{n_x} \to \mathbb{R}^{n_u}$:
\begin{align} \label{eq:dnn}
    f(x; \hat{\Theta}) = \begin{cases}
         \eta^{[\zeta]} \circ \gamma^{[\zeta]} \circ \ldots \circ \eta^{[1]} \circ \gamma^{[1]}, \,\, \zeta \geq 2\\
        \eta^{[1]} \circ \gamma^{[1]},  \,\, \zeta=1,
    \end{cases}
\end{align}
where the input of the network is $x \in \mathbb{R}^{n_x}$, the output of the network is $u \in  \mathbb{R}^{n_u}$  and $\zeta$ is the number of layers and for simplicity. If $\zeta=1$, then $f$ is describing  a \textit{shallow} neural network, else for $\zeta\geq 2$, $f$ is describing a \textit{deep} neural network. Each layer consists of a nonlinear activation function $\eta^{[k]}$ (e.g.,  \texttt{ReLU}, \texttt{sigmoid}, \texttt{HardTanh}, etc) \cite{GooBen:16} and an affine function:
$$
\gamma^{[k]}(\hat{x}^{[k-1]}) = W^{[k]} \hat{x}^{[k-1]}+ b^{[k]},
$$
where $\hat{x}^{[k-1]} \in  \mathbb{R}^{w}$ is the output of the previous layer,  \red{$w$ is the  number of neurons in each hidden layer (also referred to as the number of hidden neurons) of some arbitrary size} and $\hat{x}^{[0]} = x$.
The parameters $\hat{\Theta} = \{ \hat{\Theta}^{[1]},\ldots, \hat{\Theta}^{[\zeta]} \}$ embody   the weights and biases of the $\gamma$ function on each layer:
\begin{align*}
    \hat{\Theta}^{[k]} = \{ W^{[k]}, b^{[k]}\} \quad \forall k =\overline{1: \zeta},  
\end{align*}
where the weights are 
$$
W^{[k]} \in \begin{cases}
\mathbb{R}^{w \times n_{x}}, \quad \text{if} \,\, k=1\\
\mathbb{R}^{w \times w}, \quad \text{if} \,\, k=\overline{2: \zeta-1}\\
\mathbb{R}^{n_u \times w}, \quad \text{if} \,\, k=\zeta\\
\end{cases}
$$
and the biases are 
$$
b^{[k]} \in \begin{cases}
\mathbb{R}^{w}, \quad \text{if} \,\, k = \overline{1: \zeta-1} \\
\mathbb{R}^{n_u}, \quad \text{if} \,\, k = \zeta.
\end{cases}
$$
Once defined the neural network and selecting the training set, the goal of the deep neural network is to learn this function $f(x_{0}; \hat{\Theta}),$  that depends on the parameters $\hat{\Theta}$, such that a loss function $F(\hat{\Theta})$ (e.g. hinge loss, mean square error loss, cross entropy loss, etc.) is minimized:
$$
{\hat \Theta}^* \in \underset{\hat \Theta}{\operatorname{argmin}} \, F(\hat \Theta).$$

\noindent We define a \text{ HT deep  neural network} (\texttt{HTNN}) a feed-forward  network as described in \eqref{eq:dnn} for which all activation functions $\eta^{[k]}$ are the \texttt{HardTanh} function, which has the expression:
\begin{align}
\eta(\gamma^{[k]}, \underline{\gamma},  \overline{\gamma}) \!=\! \min \!\left(\! \overline{\gamma}, \max \!\left(\!  \gamma^{[k]}, \underline{\gamma}\right) \!\! \right) \!\!=\!\! \begin{cases}
  \underline{\gamma},\quad  \gamma^{[k]} \leq \underline{\gamma} \\
 \gamma^{[k]},  \quad \underline{\gamma}< \gamma^{[k]} < \overline{\gamma} \\
 \overline{\gamma}, \quad \gamma^{[k]} \geq \overline{\gamma}.
\end{cases}
\end{align}

\noindent The next lemma shows that  \texttt{HTNNs} yield min-max functions. 
\begin{lemma} Any neural network $f(x; \hat{\Theta})$  defined as in \eqref{eq:dnn} with  input $x\in \mathbb{R}^{n_x}$ and \texttt{HardTanh} activation functions  represents a min-max affine function. 
\end{lemma}
\begin{pf}
The \texttt{HTNN} $f(x; \hat{\Theta})$ is a min-max affine function because it only contains compositions of affine transformations with min-max affine activation functions (HardTanhs) and min-max functions are closed under addition, multiplication and composition of such functions. \qed 
\end{pf}

\noindent Next, we provide explicit bounds for the structure of a \texttt{HTNN} to be able to exactly represent a min-max affine function.

\begin{lemma} \label{lm:3}
A convex  max affine function $\mathcal{M}_{ax}: \mathbb{X} \to [ \underline{u}, \overline{u}]$,  where $\mathbb{X}$ is a compact set and $\underline{u}, \overline{u} \in \mathbb{R}$, defined as the point-wise maximum of  $\bar m$  affine functions 
$
\mathcal{A}_j(x) = c_j^T x + d_j, \,\, j={1:\bar m},$ i.e., $$ \mathcal{M}_{ax}(x)=\max_{j=1:\bar m} \mathcal{A}_j(x) = \max_{j=1:\bar m} c_j^T x+d_j,
$$
can be exactly represented by a \texttt{HTNN} with  
$$r(\bar m) = \begin{cases}
    & 0 \quad  \bar m =1\\ 
    & \bar m + r(\frac{\bar m}{2}) \quad \text{if } \bar m \text{ is even}\\
    & \bar m + r(\frac{\bar m+1}{2}) \quad  \text{if } \bar m \text{ is odd and } \bar m \neq 1, \\
\end{cases}$$ 
number of hidden neurons, $\zeta(\bar m) = \lceil \log_2 \bar m \rceil + 1$ number of layers and the maximum width $w(\bar m)=\mathbb{I}[\bar m>1] \bar m$. \end{lemma}
\begin{pf}
We start by defining the first layer of the \texttt{HTNN} as $W^{[1]} = \left[  \begin{array}{c}
     c_1 \\
     \cdots\\
     c_{\bar{m}}
\end{array}\right]\in \mathbb{R}^{\bar{m} \times n_x}$, $b^{[1]} = \left[ \begin{array}{c}
    d_1 \\
    \cdots\\
    d_{\bar{m}} 
\end{array}\right] \in \mathbb{R}^{\bar{m}}$ and $\eta^{[1]} = I_{\bar{m}}$. Thus, we obtain $ \bar{m}$ scalars $u_i \in \left[ \underline{u}_i, \, \overline{u}_i\right]$, where $ \underline{u}_i = \min_{x \in \mathbb{X}} c_i^T x + d_i$ and $\overline{u}_i = \max_{x \in \mathbb{X}} c_i^T x + d_i$ for all $i=1:\bar{m}$. In other words, we absorb the affine functions into the first layer of the HT network. Now, we need to prove for the remaining function 
$$
\gamma (u)= \max_{i =1:\bar m } u_i, \;\; \forall   u_i \in \left[ \underline{u}_i, \, \overline{u}_i \right].
$$
Our proof is based on the following  observation: 
\begin{align*}
&\max \left(u_1, u_2\right) \\
& =\max \left(0, \min \left(u_2-u_1, \overline{u}_2 - \underline{u}_1 \right)\right) + \max \left(\underline{u}_1, \min\left( u_1, \overline{u}_1 \right)\right) \\
&= \eta\left( u_2 - u_1, 0, \overline{u}_2 - \underline{u}_1\right) + \eta \left( u_1,  \underline{u}_1, \overline{u}_1\right),
\end{align*}
for any $u_1  \in \left[ \underline{u}_1, \, \overline{u}_1 \right] $ and $ u_2  \in \left[ \underline{u}_2, \, \overline{u}_2 \right]$. Indeed, the first term returns either $u_2 -u_1$ or $0$, while the second term returns always $u_1$. Using this relation recursively and similar arguments as in \cite{chen:22}, for any  $p\geq 3$ we have that:  
\begin{equation} \label{eq:maxx}
\max _{j= 1:p} u_j \!=\! \begin{cases}\max _{j=1: \frac{p}{2}} \max _{i= \{2j-1,2 j\}} u_i, &  \!\!\!\!  \text{if } p \text { is even} \\ \max _{j= 1: \frac{p+1}{2}} \phi\left(j ; u_1, u_2, \cdots, u_p\right), & \!\!\!\! \text{if } p \text{ is odd}\end{cases}
\end{equation}
for $u_j \in \left[ \underline{u}_j, \, \overline{u}_j \right], j =1:p$, where
$$
\phi\left(j ; u_1, u_2, \cdots, u_p\right)= \begin{cases}\max _{i =\{2 j-1:2 j\}} u_i, & \text{if } j = 1:\frac{p-1}{2} \\ \eta\left( u_p, \underline{u}_p, \overline{u}_p\right), & \text{if } j=\frac{p+1}{2}.\end{cases}
$$
Further, we denote by $r(p)$ the number of min-max operations from  \texttt{HardTanh} function, i.e., the number of ${\eta}(\cdot)$ evaluations, for computing the maximum of $p$ real numbers using \eqref{eq:maxx}. By expanding the operation in the aforementioned equation, one can find $r(2) = 2$ and $r(3) = 5$. We define $r(1) = 0$, since no maximum operation is needed for finding the extremum of a singleton. Thus, according to \eqref{eq:maxx}, for any integer $p \geq 3 $, we obtain the recursion:
\begin{align} \label{eq:02}
    r(p) = \begin{cases}
        p+ r(\frac{p}{2}),& \text{if } p \text{ is even}\\
        p + r(\frac{p+1}{2}), & \text{if } p \text{ is odd}.
    \end{cases}
\end{align}
Note that $r(p)$ is the number of \texttt{HardTanh} functions in a \texttt{HTNN}   that computes the maximum of $p$ real numbers or a max-affine convex function. The number of \texttt{HardTanh} functions is equivalent to the number of hidden neurons. Note that the number of \texttt{HardTanh} layers is equivalent to the number of hidden layers.  We now consider computing the maximum of $\bar m=2^p$ real numbers for any positive integer $p$. Then, every time the recursion goes to the next level in \eqref{eq:02}, the number of variables considered for computing the maximum is halved. Hence, the number of HT layers is $p$. When $\bar m$ is not a power of two, i.e., $2^p<\bar m<2^{p+1}$, then we can always construct a HT network with $p+1$ HT layers and $2^{p+1}$ input neurons, and set weights connected to the $2^{p+1}-\bar m$ "phantom input neurons" to zeros. Because $\left\lceil\log _2 \bar m\right\rceil=p+1$ for $2^p<\bar m<2^{p+1}$, the number of HT layers is $\left\lceil\log _2 \bar m\right\rceil$ for any positive integer $\bar m$. Hence, we have  $\zeta(\bar m)=\left\lceil\log _2 \bar m\right\rceil+1$. Note that $r(\bar m)$ is a strictly increasing sequence. Therefore, the maximum width of the network is given by the width of the first hidden layer. When we have one layer  or $\bar m=1$, the width is 0. When the number of layers is greater than 1 or $\bar m>1$, we have $w(\bar m) = \bar m$. \qed 
\end{pf}
For the reader's convenience, we introduce 2 lemmas regarding the neural networks connections from \cite{chen:22}, that we will use in our final result (for proof see Lemma 4 and Lemma 8 from \cite{chen:22}, respectively).  
\begin{lemma} \label{lm:2nn}
    Consider two feed-forward neural networks $f^{[1]}$ and $f^{[2]}$. Each network is characterized by $r^{[i]}$ number of hidden neurons, $w^{[i]}$ maximum width and $\zeta^{[i]}$ number of layers, with $i=1:2$. Then, there exist a feed-forward neural network, denoted $f$, that represents the composition of two neural networks, i.e., $f = f^{[1]} \circ f^{[2]}$, with $\zeta = \zeta^{[1]} + \zeta^{[2]} - 1$, $r =r^{[1]} + r^{[2]} $ and $w = \max(w^{[1]}, w^{[2]})$.
\end{lemma}
\begin{lemma}\label{lm:l-nn}
    Consider $l$ feed-forward neural networks $f^{[1]}, \ldots, f^{[l]}$, each network being characterized by $r^{[i]}$ number of hidden neurons, $w^{[i]}$ maximum width, $\zeta^{[i]}$ number of layers and $m_i$ the dimension of the NN output, with $i=1:l$. Then, there exist a feed-forward neural network, denoted $f$, that represents the parallel connection of the $l$ neural networks, i.e., $f(x) = \left[\begin{array}{c}
    f^{[1]}(x)\\
    \ldots \\
    f^{[l]}(x) 
    \end{array}\right],$ with:
    \begin{align*}
        & \zeta = \max_{i=1:l} \zeta^{[i]} , \qquad
        w = \sum_{i=1}^{l} \max(w^{[i]}, 2m_i)\\
        & r = \sum_{i=1}^{l} r^{[i]} + 2m_i( \zeta - \zeta^{[i]}).
    \end{align*}
\end{lemma}
Now, based on the previous lemmas, we prove that any scalar min-max function can be exactly represented by a \texttt{HTNN}.
\begin{theorem}\label{th:01}
Any min-max function $\mathcal{M}: \mathbb{X} \to \left[ \underline{u}, \overline{u}\right]$, with $\mathbb{X}$ a compact set and $\underline{u}, \overline{u} \in \mathbb{R}$, defined by  $m$  affine functions $\mathcal{M}_j(x) = c_j^T x + d_j, \,\, j={1:m}$, i.e., 
 \begin{align}
\label{eq1}
 \mathcal{M}(x) = \min_{1 \leq i \leq l}  \max_{j \in \mathcal{C}_i}   \mathcal{M}_j(x) = \min_{1 \leq i \leq l}  \max_{j \in \mathcal{C}_i} c_j^T x+d_j
 \end{align}
  can be exactly represented by a \texttt{HTNN} whose number of layers $\zeta$, maximum width $w$ and number of hidden neurons $r$ satisfy:
  \begin{align*}
 & \zeta \leq \lceil \log_2 l \rceil +  \lceil \log_2 m \rceil + 1 \\
 & w \leq \mathbb{I}[l>1] l \max(m, 2) \\
 & r \leq 2l \left( 1+ 2m + \lceil \log_2 m\rceil \right) -2.
  \end{align*}
\end{theorem}
\begin{pf}
According to \eqref{eq1} there are $l$ minima  to be computed which are the maximum of $|\mathcal{C}_i|$ real numbers with $i =1:l$. Then, the value of $\mathcal{U}$ can be computed by taking the minimum of the resulting $l$ maxima. Below we prove that these operations can be achieved by a \texttt{HTNN}. From Lemma \ref{lm:3} one can compute the maximum of $\bar m$ numbers given by $\bar m$ affine functions with a HT network with $r(\bar m)$ hidden neurons, $\zeta(\bar m)$ layers  and $w(\bar m)$ maximum width. We construct the HT network for $\mathcal{M}(x)$ as follows: we create $l$ \texttt{HTNNs} and arrange them in parallel each one corresponding to a convex max affine function $\max_{j  \in \mathcal{C}_i} \mathcal{M}_j(x)$ for $i=1:l$, using Lemma \ref{lm:3}. If we denote  $|\mathcal{C}_i| = \bar m_i$, then each individual network has $r(\bar m_i)$ hidden neurons, $\zeta(\bar m_i)$ layers  and $w(\bar m_i)$ maximum width, for all $i=1:l$. These $l$ concatenated \texttt{HTNNs} have the same input $x$ and outputs $l$ real numbers. Finally, we create a \texttt{HTNN} with $r(l)$ hidden neurons, $\zeta(l)$ layers and $w(l)$ maximum width, that takes the minimum of $l$ real numbers using again Lemma \ref{lm:3} since $\min(\cdot) = -\max(-\cdot)$. Combining Lemma \ref{lm:2nn} and \ref{lm:l-nn}, the number of layers of this extended \texttt{HTNN} is:
$$\zeta = \zeta(l) + \max_{i=1:l} \zeta(\bar m_i) -1 .$$
The last term appears due to the fact that the first layer  $\zeta(l)$ of HTNN is absorbed  
into the last layer of the other one, since their dimensions are compatible.
Further, the number of hidden neurons is 
\begin{align}\label{eq:nhn}
    r = r(l)+ \sum_{i =1}^{l} r( \bar{m}_i) + 2 \left(\max_{j=1:l} \zeta(\bar{m}_j) - \zeta(\bar{m}_i)\right),
\end{align}
and the maximum width is
\begin{align*}
   w = \max\left(\sum_{i=1}^{l} \max(w(\bar m_i), 2), w(l) \right).
\end{align*}
Since $ \mathcal{C}_i \subseteq \{ 1, \ldots, m \} $ it follows that  $| \mathcal{C}_i| \geq 1$ and $| \mathcal{C}_i| = \bar{m}_i \leq m$ for all $i =1:l$.
Therefore the number of layers can be  bounded from above by (see Lemma \ref{lm:3})
$$ \lceil \log_2 l \rceil +  \lceil \log_2 m \rceil + 1.$$
Similarly we can derive an upper bound for the maximum width: 
$$
\max( l \max (w(m), 2), w(l) ) \leq \max( l\max(m,2), l) = l \max(m, 2). 
$$
For the hidden number of neurons, let us first bound the $r(p)$ sequence from \eqref{eq:02}. Since $r(p)$ is  strictly increasing, then:
\begin{align}\label{eq:rbound}
r(p) & = r(2^{\log_2p}) \leq r(2^{\lceil \log_2p \rceil}) \overset{\eqref{eq:02}}{=} \sum_{i=1}^{\lceil \log_2p \rceil} 2^i 
 = 2^{\lceil \log_2p \rceil +1} -2 \nonumber \\ 
 & = 2 (2^{\lceil \log_2p \rceil} -1) < 2 (2^{ \log_2p +1} -1) = 2(2p-1).
\end{align}
Using $| \mathcal{C}_i| \geq 1$ and $| \mathcal{C}_i| = \bar{m}_i \leq m$ in \eqref{eq:nhn}, we obtain:
\begin{align*}
    r(l) + \sum_{i =1}^{l} r( m) + 2 \left( \zeta(m) - \zeta(1)\right) =
   r(l) + l(r(m) + 2\lceil \log_2 m\rceil)\\
   \leq 2l \left( 1+ 2m + \lceil \log_2 m\rceil \right) -2,
\end{align*}
where the last inequality is given by \eqref{eq:rbound}. This proves our statements. \qed
\end{pf}
Note that in the previous theorem, we treat the scalar case. Thus, we provide next the extension to the vector case.
\begin{theorem} \label{th:02}
Any min-max function $\mathcal{M}: \mathbb{X} \to [\underline{\ub}, \, \overline{\ub}]$, with $\underline{\ub}, \, \overline{\ub}\in \mathbb{R}^{n_u}$ and $\mathbb{X}$ a compact set, defined as:
\begin{align*}
\label{eq2}
 \mathcal{M}(x) = \min_{1 \leq i \leq l}  \max_{j \in \mathcal{C}_i} C_j x+d_j,
 \end{align*}
with $\mathcal{C}_1, \ldots, \mathcal{C}_l \subseteq\{1, \ldots, m\}$, can be exactly represented by a parallel connection of $n_u$ \texttt{HTNNs}:
\begin{equation}
   \left[\begin{array}{cc}
       f^{[1]}(x; \hat{\Theta}^{[1]})  \\
       \vdots \\
        f^{[n_u]}(x; \hat{\Theta}^{[n_u]}) 
   \end{array} \right],
   \label{eq:vectorNN}
\end{equation}
where each \texttt{HTNN}, $f^{[i]}$, has $\zeta_i$ number of layers, $r_i$ number of hidden neurons and $w_i$ maximum width, for $i = 1:n_u$. Then, the $n_u$ parallel \texttt{HTNN} has:
\begin{align*}
        & \zeta = \max_{i=1:n_u} \zeta_i , \qquad
        w = \sum_{i=1}^{n_u} \max(w_i, 2),\\
        & r = \sum_{i=1}^{n_u} r_i + 2( \zeta - \zeta_i).
    \end{align*}
\end{theorem}
\begin{pf}
    A min-max function $\mathcal{M}: \mathbb{X} \to [\underline{\ub}, \, \overline{\ub}]$ can be divided  into one scalar min-max function per output dimension:
\begin{align*}
 \mathcal{M}_k(x): \mathbb{X} \to [\underline{u}, \, \overline{u}] \;\;\; \forall  k =1:n_u.
\end{align*}
Applying Theorem \ref{th:01}, each scalar min-max function can be exactly represented by a HTNN, denoted $ f^{[k]}(x; \hat{\Theta}^{[k]}), \forall k=1:n_u$. Using a parallel architecture as in \eqref{eq:vectorNN}, the min-max function $\mathcal{M}$ is vectorized. Using Lemma \ref{lm:l-nn} we obtain the number of layers, maximum width and the number of hidden neurons of the network as stated in the theorem.  \qed
\end{pf}
From Theorems \ref{th:01} and \ref{th:02} it follows that we can represent exactly a $n_u$ vector min-max function $ \mathcal{M}(x) = \min_{1 \leq i \leq l}  \max_{j \in \mathcal{C}_i} C_j x+d_j,$ with $\mathcal{C}_1, \ldots, \mathcal{C}_l \subseteq\{1, \ldots, m\}$, by a \texttt{HTNN} with: 
\begin{align*}
    &\text{number of layers } \zeta \leq \mathcal{O}(\log_2 l + \log_2 m)\\ 
    &\text{maximum width } w \leq \mathcal{O}(n_u l m) \\
    & \text{number of hidden neurons } r \leq \mathcal{O}(n_u l m).
\end{align*}
\noindent Hence, the complexity of the network  is linear in the dimension of the output of the network ($n_u$),  in the number of affine terms ($m$) and also  in the number of min terms ($l$). In particular,  this shows that the number of hidden neurons required to exactly represent a min-max function is at most an expression given by the number of pieces times number of min terms. 

\begin{remark}
{ Recently, in \cite{chen:22} it was proved that the number of hidden neurons required to exactly represent a scalar PWA function with a deep neural network using  ReLU is a quadratic function of the number of pieces ($m$). On the other hand, it was shown in \cite{XuBoo:16} that in general the min-max representation is more efficient than the corresponding PWA representation, i.e., $l \ll m$. Hence,  when $l \leq m$ our approach of learning  a PWA/min-max function through an HTNN  network is more beneficial (e.g., in terms of number of hidden neurons) than learning through a  ReLU type  neural network.} 
\end{remark}

\noindent Further, our previous results can be used to learn exactly an MPC law, considering  a training set  $\mathcal{S} \!=\! \left\{ \left(\ub^{*}_s, x_{0,s}\right) \!\mid \!s \!=\!1: D\right\}$, where $\ub^{*}_s = \ub^{*}(x_{0,s})$ is the optimal command and $x_{0,s} \in \mathbb{X}$ is the associated initial state of the system. Then, the goal of the \texttt{HTNN} is to learn exactly a min-max function  $\ub^{ \texttt{\scriptsize HTNN}}(x_0)$  depending on  the parameters $\hat{\Theta}$ such that the following mean-squared error  empirical loss is minimized:
$$
\hat{\Theta}^* \in \underset{\hat \Theta}{\operatorname{argmin}} \, F(\hat \Theta):=\frac{1}{D} \sum_{s=1}^D\left\|\ub^{*}_s - \ub^{ \texttt{\scriptsize HTNN}}(x_{0,s}; {\hat \Theta}) \right\|^2. 
$$
Given the structure of the \texttt{HTNN} network,  we obtain again a min-max representation for the learned function $\ub^{ \texttt{\scriptsize HTNN}}(x_0)$, which depends on the (optimal) parameters  $\hat{\Theta}^*$ of the \texttt{HTNN}: 
\begin{align}
\label{mm-htnn}
\mathbf{u}^{\texttt{\scriptsize HTNN}}(x_0) &= \mathcal{M}_{in}^{\texttt{\scriptsize HTNN}}\mathcal{M}_{ax}^{\texttt{\scriptsize HTNN}}(x_0; \hat{\Theta}^*).    
\end{align} 
From previous theorems it follows that $\ub^{ \texttt{\scriptsize HTNN}}(x_0) = \ub^*(x_0)$, provided that the number of layers and hidden neurons are chosen as in Theorem 2.


\section{Unfolded HardTanh deep neural networks (U-HTNN) for approximate representation of  MPC laws}
\noindent In this section, in order to eliminate the disadvantages of generic  black-box  \texttt{HardTanh} deep neural networks, we  propose unrolled deep learning techniques, which bring the physics of the model and standard optimization techniques into the architecture design. Hence, we consider approximating the MPC law using \texttt{HardTanh} based deep neural networks whose architecture is  inspired from first order optimization methods, we call this strategy deep unfolding first order methods-based architectures. 
Note that the derivations from below are similar to \cite{lupnec:23}. However section 4.1 is new compared to \cite{lupnec:23}. Similar to unfolded denoiser for images in \cite{LePus:22}, we design the architectures of our MPC deep neural networks, taking  inspiration from the first order methods, PGD and APGD. Since PGD can be recovered from APGD by choosing $\beta_k=0$, we will consider in the sequel only  APGD algorithm. First, let us replace the expression of $\yb_{k}$  in the update step of $\ub_{k+1}$:
\begin{equation*}
    \ub_{k+1} \!=\! \Pi_{\mathbb{U}^N}\left( (\operatorname{I_{\nu}} - \alpha_k \mathbf{H}) \left[ (1+\beta_{k-1}) \ub_k \!- \beta_{k-1} \ub_{k-1}\right] - \alpha_k  \mathbf{q}(x_0)\right).
\end{equation*}
Observing that $\ub_{k+1}$ depends on the current and previous iterates, we rearrange the relation in a matrix format:
\begin{align} \label{eq:1}
&\left[\begin{array}{c}
\ub_{k}  \\
\ub_{k+1} \\
\end{array}\right] =  \left[\begin{array}{c}
0  \\
- \alpha_k \mathbf{q}(x_0) \\
\end{array}\right] +  \nonumber\\
&\left[\begin{array}{cc}
0 & \operatorname{I_{\nu}}\\
-\beta_{k-1}( \operatorname{I_{\nu}}- \alpha_k \mathbf{H}) & (1+\beta_{k-1}) (\operatorname{I_{\nu}} -\alpha_k \mathbf{H}) \\
\end{array}\right] \left[\begin{array}{c}
\ub_{k-1}  \\
\ub_{k} \\
\end{array}\right].
\end{align}
Further, the projection operator is applied only on $\ub_{k+1}$. Setting $\ub_0 = \ub_1$ results in $\yb = \ub_1$ and $\ub_2 = (\operatorname{I}_{\nu} - \alpha_k \mathbf{H}) - \alpha_k \mathbf{q}(x_0)$. Due to this setting, the parameters of the first layer will have a different structure than the rest of the network.  We design the rest of the layers of our neural network architecture using the expression \eqref{eq:1}. Since the output of the layer, contains $\ub_{k}$ and $\ub_{k+1}$, the last layer extracts $\ub^{ \texttt{\scriptsize U-HTNN}}(x_0;{\hat \Theta}) = \ub_{\zeta-1}$, where $\zeta$ is the number of layers. Recalling that  $\mathbf{q}(x_0) = \mathbf{A}^T \mathbf{Q} \mathbf{A}_N x_0$, we denote with $Q_1 = \mathbf{H}$ and $Q_2 = \mathbf{A}^T \mathbf{Q} \mathbf{A}_N$. \red{Thus, the weights, biases and activation functions for the layers of our deep unfolded network, called \texttt{U-HTNN}, are defined, taken into account the structure of the MPC objective function and constraints, as  follows}:
\begin{align*}
    & W^{[1]}=\left[\begin{array}{c}
\operatorname{I_{\nu}} \\
\operatorname{I_{\nu}} -\alpha^{[1]} Q_1^{[1]}
\end{array}\right], \,\,\\
& b^{[1]}=\left[\begin{array}{c}
0 \\
- \alpha^{[1]}Q_2^{[1]} x_{0}
\end{array}\right], \,\eta^{[1]}=\left\{ \!\begin{array}{c}
\mathcal{I}_{d} \\
\bar{\eta}(\cdot \;, \overline{\ub}, \underbar{$\ub$})
\end{array} \! \right\}\,, 
\end{align*}
\begin{align*}
& W^{[k]}=\left[\begin{array}{cc}
0 & \operatorname{I_{\nu}} \\
- \beta_{k-1} \left(\operatorname{I_{\nu}} - \alpha^{[k]}Q_1^{[k]} \right) & (1+ \beta_{k-1}) \left(\operatorname{I}_{\nu} - \alpha^{[k]}Q_1^{[k]}\right)
\end{array}\right]\\
& b^{[k]} \!=\! \left[\!\!\begin{array}{c}
0 \\
- \alpha^{[k]}Q_2^{[k]} x_{0}
\end{array}\!\!\right],\eta^{[k]}\!=\!\left\{ \!\!\! \begin{array}{c}
\mathcal{I}_{d} \\
\bar{\eta}(\cdot \;, \overline{\ub}, \underbar{$\ub$})
\end{array} \!\!\! \right\}\, \forall k \!\in\! \{2, \ldots, \zeta\!-\!1\}, \\ 
& W^{[\zeta]}= \operatorname{I_{2\nu}}, \;\;   b^{[\zeta]}= 0, \;\;  \eta^{[\zeta]}=\left\{ \!\begin{array}{c}
0 \\
\mathcal{I}_{d}
\end{array} \! \right\}\,,
\end{align*}

\noindent where  $\mathcal{I}_d$ denotes the identity function. Observe that the first components in the parameters and activation functions, are chosen in a manner to propagate $\ub_k$. Recall that $\eta (\cdot)$ denotes the standard activation function \texttt{HardTanh}. Note that in this network we need to learn the entries of the matrices ${Q}_1^{[k]}$ and $ {Q}_2^{[k]}$ and also $\alpha^{[k]}$, $\beta_k$, for $k=1:\zeta-1$. The extra parameter $\beta_k$ allows us to consider the basic optimization  algorithm PGD ($\beta_k =0$) or its accelerated version APGD.   After training our  \texttt{U-HTNN}, we get an approximation map $\ub^{ \texttt{\scriptsize U-HTNN}}(x_0)$, \text{with}   ${\hat \Theta} = (\hat{Q}_1^{[k]}, \hat{Q}_2^{[k]})_{k=1:\zeta}$.  Given the structure of the \texttt{U-HTNN} network,  we again obtain  a min-max representation for the learned function $\ub^{ \texttt{\scriptsize U-HTNN}}(x_0)$, which depends on the (optimal) parameters  $\hat{\Theta}^*$ of the \texttt{U-HTNN}: 
\begin{align}
\label{mm-uhtnn}
\mathbf{u}^{\texttt{\scriptsize U-HTNN}}(x_0) &= \mathcal{M}_{in}^{\texttt{\scriptsize U-HTNN}}\mathcal{M}_{ax}^{\texttt{\scriptsize U-HTNN}}(x_0; \hat{\Theta}^*).    
\end{align}
Denote $\zeta$ number of layers and $w$ the number of neurons in
each hidden layer. From existing  universal approximation theorems for PWA/min-max functions one can also derive a uniform bound on this approximation, see  \cite{KarLuc:20,HorSti:89} for more details: 
\begin{align*}
\| \mathbf{u}^{*}(x_0) -  \ub^{ \texttt{\scriptsize U-HTNN}}(x_0) \| & \leq \mathcal{O} \left(\zeta,w \right) \;\;\;  \forall x_0 \in \mathbb{X}.
\end{align*} 

\noindent The following proposition establishes the relation between \texttt{U-HTNN}  and (A)PGD algorithms.  The proof follows from the fact that the projection operator fits the {\tt HardTanh} activation function.

\begin{theorem}
Set, for every $k \in\{1, \ldots, \zeta\}, Q_1^{[k]} \in \mathbb{R}^{\nu\times \nu}$,   $Q_2^{[k]} \in \mathbb{R}^{\nu \times \nu}, \;  W^{[k]}, \; b^{[k]}$ and $\eta^{[k]}$ provided by  \texttt{U-HTNN}. If $ Q_1^{[k]}= \alpha_k \mathbf{H}$, \;  $Q_2^{[k]}= \mathbf{A}^T \mathbf{Q} \mathbf{A}_N$ and $\ub_0= \ub_1= \yb_1$, then  \texttt{U-HTNN}  fits the generic (A)PGD  optimization schemes \eqref{pgd} or \eqref{apgd}.
\end{theorem}


\subsection{Structured \texttt{U-HTNN} architecture}
\label{sec:suhtnn}
\noindent In this subsection, we impose more structure on the weights of the \texttt{U-HTNN} by considering the particular form  of $\mathbf{H} = \mathbf{R} + \mathbf{A}^T\mathbf{Q}\mathbf{A}$ and $\mathbf{q}(x_0) = \mathbf{A}^T\mathbf{Q}\mathbf{A}_{N} x_0$ in the MPC optimization problem  \eqref{eq:qp}. For reducing the number of parameters, we consider the case $\mathbf{R} = I$, which is not a restrictive assumption. 
Notice that in both $\mathbf{H}$ and $\mathbf{q}(x_0)$  the term $\mathbf{A}^T\mathbf{Q}$ appears, which we will denote with $Q_{11}^{[k]}$. The remaining terms are denoted  $Q_{12}^{[k]} = \mathbf{A}$ and $Q_{21}^{[k]} = \mathbf{A}_N$. With these notations the structured \texttt{U-HTNN} (called \texttt{S-U-HTNN}) is  parameterized as follows:
\begin{align*}
    & W^{[1]}=\left[\begin{array}{c}
\operatorname{I}_\nu \\
\operatorname{I}_\nu -\alpha^{[1]} Q_{11}^{[1]}Q_{12}^{[1]}
\end{array}\right], \, b^{[1]}=\left[\begin{array}{c}
0 \\
- \alpha^{[1]}Q_{11}^{[1]}Q_{21}^{[1]} x_{0}
\end{array}\right], \\
& W^{[k]}=\left[\begin{array}{cc}
0 & \operatorname{I}_\nu \\
W_1^{[k]}& W_2^{[k]}
\end{array}\right], \; b^{[k]} = \left[\begin{array}{c}
0 \\
- \alpha^{[k]}Q_{11}^{[k]}Q_{21}^{[k]} x_{0}
\end{array}\right], \\
& \forall k \!\in\! \{2, \ldots, \zeta\!-\!1\}, \; W^{[\zeta]}= \operatorname{I}_{2\nu}, \;\;   b^{[\zeta]}= 0,
\end{align*}
where $W_1^{[k]} = - \beta_{k-1} (1-\alpha^{[k]})  \operatorname{I}_\nu + \beta_{k-1} \alpha^{[k]}Q_{11}^{[k]} Q_{12}^{[k]}$ and $W_2^{[k]} = (1+\beta_{k-1}) (1-\alpha^{[k]})  \operatorname{I}_\nu - (1+\beta_{k-1}) \alpha^{[k]} Q_{11}^{[k]} Q_{12}^{[k]} $. 
The activation functions have the same structure as in \texttt{U-HTNN}. 

\medskip 

\noindent Examining even further the matrices appearing in the dense MPC problem  \eqref{eq:qp}, we can have a super structured architecture. In this version, called \texttt{SS-U-HTNN}, the terms $Q_{11}^{[k]}$ and $Q_{12}^{[k]}$ for all $k=1:\zeta-1$ are taken  quasi-upper and -lower triangular matrices, respectively, since $\mathbf{A}$ is a quasi-lower triangular matrix and $\mathbf{Q}$ is a block diagonal matrix. As for the parameter $Q_{12}^{[k]}$, given the form of $\mathbf{A}_N$ (see \eqref{eq:mpc1}), we will only learn the small matrix $A\in \mathbb{R}^{n_x \times  n_x}$. Note that the number of  parameters to be learned now in \texttt{SS-U-HTNN} is given by the number of states of the system and thus it  is much smaller than in the previous proposed networks \texttt{U-HTNN} and \texttt{S-U-HTNN}, respectively.


\subsection{Robustness of \texttt{HTNN} and \texttt{U-HTNN} architectures}
\noindent Following \cite{ComPes:20,LePus:22},  \texttt{HTNN} and \texttt{U-HTNN}  have Lipschitz
behavior with Lipschitz constant $\textbf{L}=\prod_{k=1}^\zeta\left\|W^{[k]}\right\|$. Thus,
given an input $x_0$ and some perturbation $\epsilon$, we can majorize the
perturbation on the output via the inequality:
$$
\left\|   \ub^{ \texttt{\scriptsize NN}}(x_0 + \epsilon, \hat{\Theta}) -  \ub^{ \texttt{\scriptsize NN}} (x_0, \hat{\Theta})\right\| \leq \textbf{L} \, \|\epsilon\|,
$$
where $\ub^{ \texttt{\scriptsize NN}} \in \{ \ub^{ \texttt{\scriptsize HTNN}}, \ub^{ \texttt{\scriptsize U-HTNN}}\}$. Therefore,  Lipschitz constant $\textbf{L}$ can be used as a certificate of  robustness of our networks provided that it is tightly estimated. Tighter bounds for  Lipschitz constant exist but at the price of more complex computations, see e.g., \cite{FG:2021}.


\section{Numerical experiments}

\begin{figure}[!htbp] 
\centering
	\begin{subfigure}{0.5\textwidth} 
        \centering
		\includegraphics[scale=0.7]{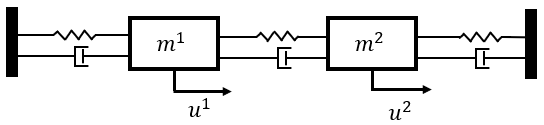}
		\caption{} 
		\label{SubFig:Jasper_original_non_sparisifed}
	\end{subfigure}	
	\begin{subfigure}{0.5\textwidth} 
 \centering
		\includegraphics[scale=0.6]{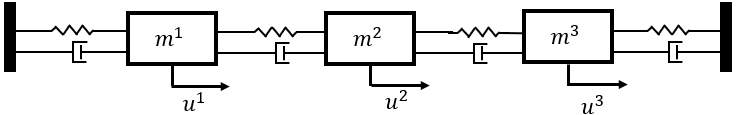}
		\caption{} 
		\label{SubFig:Jasper_sparisifed}
	\end{subfigure}	
	\caption{System of 2 and 3  oscillating masses connected each other through pairs of spring-damper blocks, and to walls (dark blocks on the sides).}
 \label{fig:sys}
\end{figure}

\begin{table*}[h]
\centering
\begin{tabular}{|c|c|c|cc|cc|cc|cc|}
\hline
\multirow{2}{*}{System} & \multirow{2}{*}{MPT} & \multirow{2}{*}{APGD} & \multicolumn{2}{c|}{HTNN}                                      & \multicolumn{2}{c|}{U-HTNN}                  & \multicolumn{2}{c|}{S-U-HTNN}                & \multicolumn{2}{c|}{SS-U-HTNN}               \\ \cline{4-11} 
                         &                      &                       & \multicolumn{1}{c|}{$\zeta=7,\, w = 4,$} & $ \zeta=10, \, w = 8$ & \multicolumn{1}{c|}{$ \zeta=3$} & $ \zeta=7$ & \multicolumn{1}{c|}{$ \zeta=3$} & $ \zeta=7$ & \multicolumn{1}{c|}{$ \zeta=3$} & $ \zeta=7$ \\ \hline
(a)                 & 0.52                & 0.09                 & \multicolumn{1}{c|}{0.03}                    & 0.07              & \multicolumn{1}{c|}{0.08}        & 0.12      & \multicolumn{1}{c|}{0.09}       & 0.17      & \multicolumn{1}{c|}{0.08}      &  0.12           \\ \hline
(b)                 & *                    & 0.12                 & \multicolumn{1}{c|}{0.06}                    &  0.1                   & \multicolumn{1}{c|}{0.11}           & 0.14            & \multicolumn{1}{c|}{0.1}           &   0.19        & \multicolumn{1}{c|}{0.1 }         &  0.15          \\ \hline
\end{tabular}
\caption{Average time (sec) for online evaluation of the proposed (approximate) MPC laws for systems of 2 and 3  oscillating masses. }
\label{tb:01}
\end{table*}
\noindent In this section, we illustrate the performance of the proposed (approximate) MPC laws learned with  general/unfolded  neural networks and compare them  with those computed with optimization methods.  The goal is to stabilize the system of oscillating masses, see \cite{FG:2021}. The masses-springs-dampers configurations we study are given in Figure \ref{fig:sys}, i.e.,  we consider a system with 4 or 6 states and 2 or 3 inputs, respectively. For numerical simulations, all masses are 1, springs constants 1, damping
constants 0.5 and we discretize the dynamics of the system using a sampling rate  0.1. For each mass $m^i$, with $i\in \{1,2,3\}$, the associated state $x^i \in \mathbb{R}^2$ contains the position and the velocity of the mass.  The control inputs are subject to the constraint $|u^i| \leq 1$. 
To complete the MPC formulation \eqref{eq:mpc}, we set the receding horizon $N = 5$ and the cost matrices $Q$ and $R$ to identity.  In order to compute the explicit MPC law, we are using the MPT toolbox \cite{mpt} with Matlab R2020. For the first order algorithm, APGD, we set $\alpha$ and $\beta$ according to \cite{Nes:04} (scheme III). For learning the (approximate) MPC laws generated by the general  neural network   \texttt{HTNN}, unfolded neural network  \texttt{U-HTNN} and unfolded neural networks with structure  \texttt{S-U-HTNN}/\texttt{SS-U-HTNN}, respectively, \red{we randomly generated a  set of $500$ initial  points  in the interval $[-4, 4] $ for the mass position and $[-10, 10]$ for the velocity, for each $i\in \{1,2,3\}$. For each initial point, the MPC problem is solved  using  APGD algorithm over a simulation  horizon of length 100, yielding an optimal sequence of MPC inputs of length $N=5$ for each simulation step. We make use of the entire optimal sequence of MPC inputs and associated with the corresponding  initial state over the simulation horizon. Thus, we obtain  50000  pairs $( x_0, \ub^*(x_0))$. We use 70\% of the data for training, 10\% for validation and 20\% for testing.} For training the neural networks we used the ADAM optimizer with a learning rate  $10^{-3}$.


\begin{figure}[h]
\includegraphics[width=9.5cm, height=4.5cm]{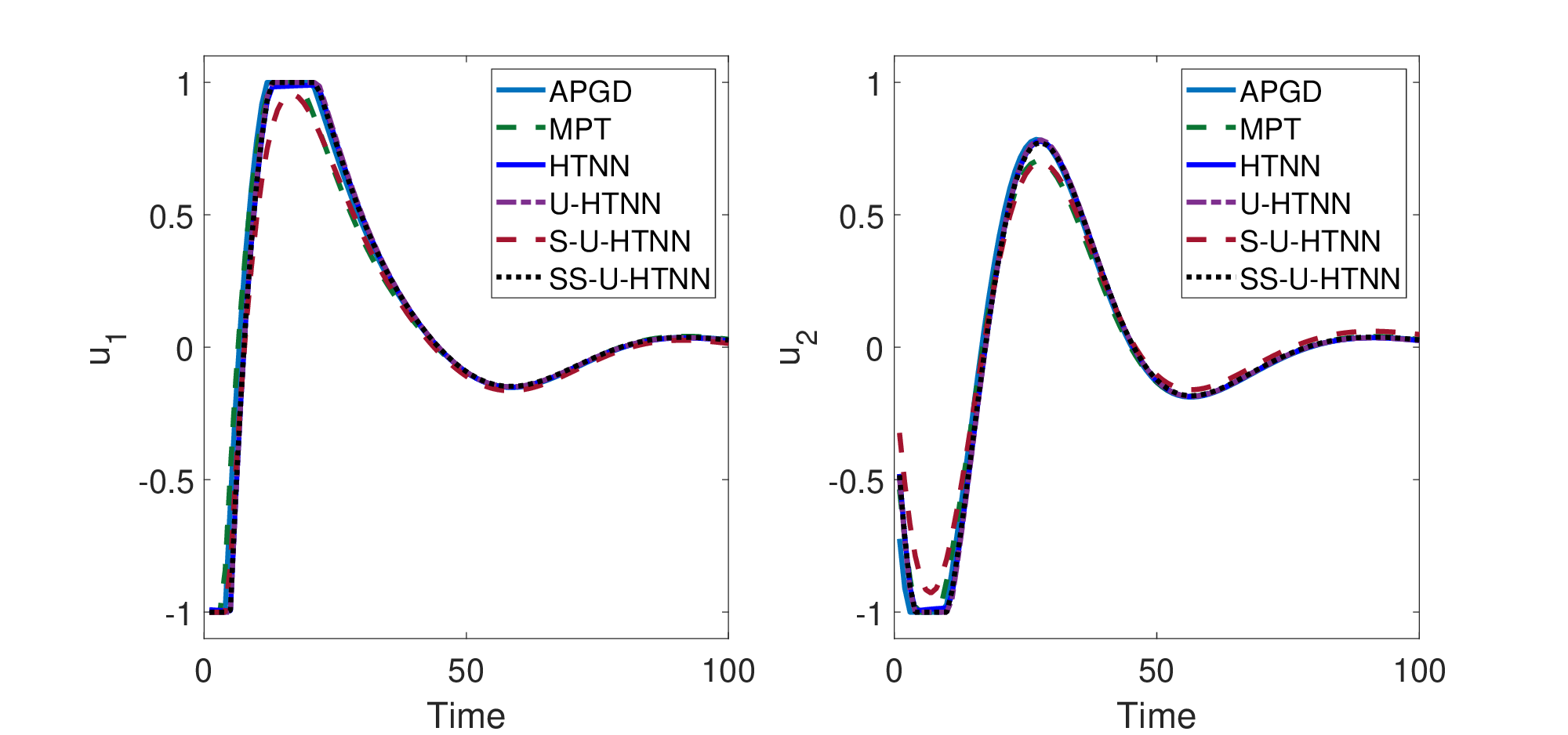}
\includegraphics[width=9.5cm, height=4.5cm]{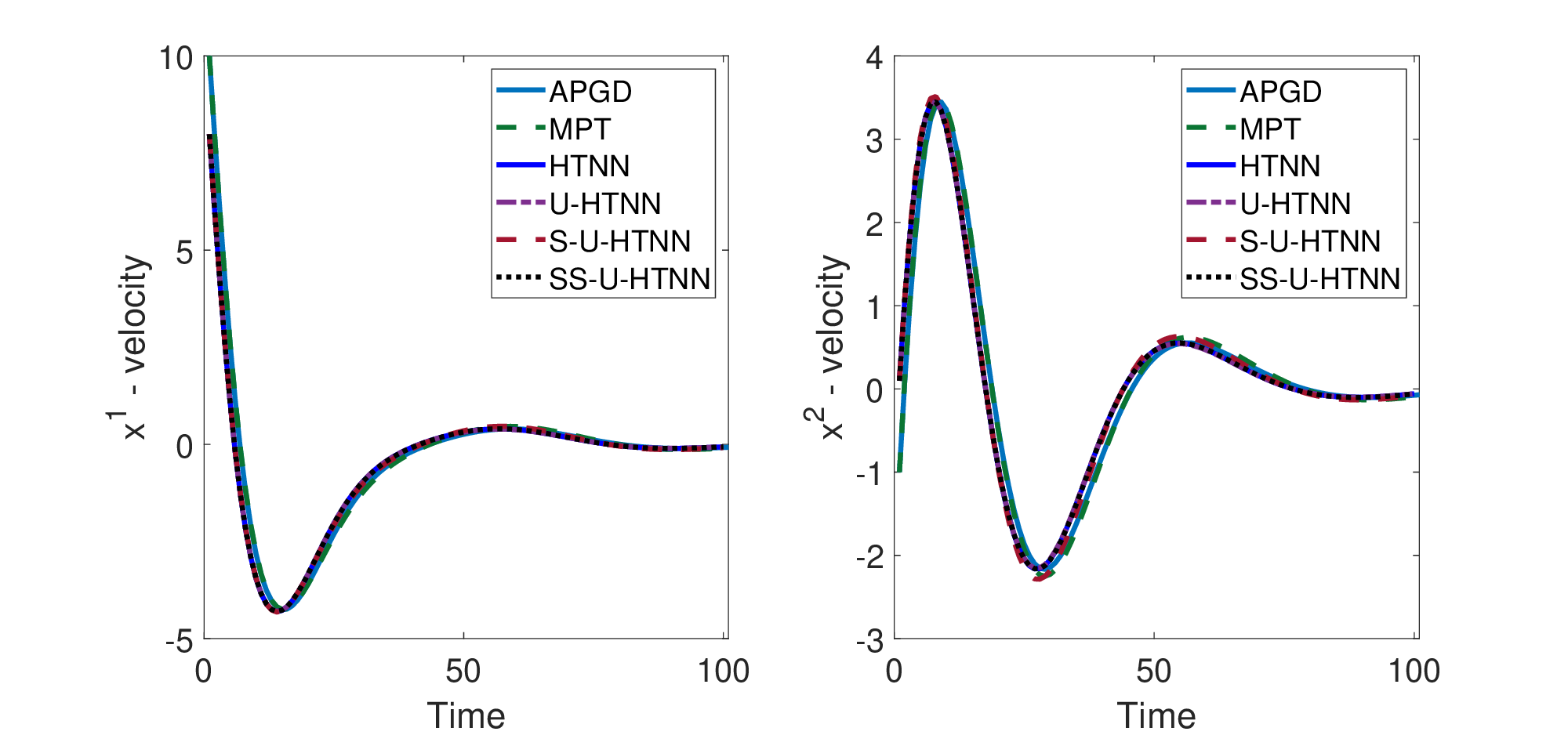}
\caption{Closed-loop  oscillating 2 masses system (a): comparison between explicit MPC (MPT), APGD and learned MPC laws using \texttt{HTNN}, \texttt{U-HTNN}, \texttt{S-U-HTNN} and \texttt{SS-U-HTNN} networks  - inputs (top) and  states (bottom) trajectories for the initial state $x_0 = [4,\, 10,\, -1,\, -1]^T$. }
\label{fig:traj}
\end{figure}

\medskip  

\noindent From Figure 2 we observe that  all MPC inputs and states  trajectories  corresponding to neural network/optimization methods  are very similar  and  keep all the controllers between the imposed bounds. Moreover, the computed MPC inputs ($u^1$ and $u^2$) are also on the boundary  for some time in all the approaches we have considered in simulations.     

\medskip 

\noindent Further, in Table \ref{tb:01}, we compare  the average time for online evaluation of all MPC laws corresponding to neural network/ optimization methods.  For the system of 2 oscillated masses (a), MPT toolbox yields 4647 polyhedral partitions, while for the system of 3 oscillated masses (b), the MPT toolbox was not able to find the control law in less than one hour. For the neural network approaches we consider different sizes for the number of layers ($\zeta$) and the maximum width ($w$). Note that for online evaluation in APGD and general \texttt{HTNN} we perform only matrix-vector operations, while in (structured) unfolded \texttt{HTNN} networks we need additionally matrix-matrix computations. Hence, for (structured) unfolded \texttt{HTNN} networks the online time may be longer, however the architecture (e.g., number of layers) can be much simpler than for the general \texttt{HTNN}. E.g., as one can see from Figure \ref{fig:traj}, \texttt{U-HTNN} network with $\zeta=3$ layers performs as good as the general \texttt{HTNN}  with $\zeta = 7$ layers. However, from Table \ref{tb:01} one can see that the neural network approaches are much faster than the classical multi-parametric methods (such as MPT) and slightly faster than the online accelerated gradient methods (such as APGD).



\section{Conclusions}
\noindent In this paper we have designed  \texttt{HardTanh} deep  neural network architectures for exactly or approximately computing  linear MPC laws. We first considered generic black-box networks and derived theoretical bounds on the minimum number of hidden layers and neurons such that it can exactly represent a given MPC law. We have also designed unfolded  \texttt{HardTanh} deep  neural network architectures based on first order  optimization algorithms to approximately compute  MPC laws.  Our \texttt{HardTanh} deep neural networks have a  strong potential to be used as approximation methods for MPC laws, since they usually yield better approximation  and  less time/memory requirements than previous approaches, as we have illustrated via simulations.


\section*{ACKNOWLEDGMENT}
\noindent The research leading to these results has received funding from UEFISCDI PN-III-P4-PCE-2021-0720, under  project L2O-MOC, nr. 70/2022.


\end{document}